\documentclass{article}

\usepackage{Jared-MK}
\usepackage{graphicx} 
\usepackage{mathtools} % http://ctan.org/pkg/mathtools
\usepackage{hyperref}
\usepackage{thmtools}
\usepackage{thm-restate}
\usepackage{cleveref}

\newcommand{\n}{\nabla}

\newtheorem*{theorem*}{Theorem}

\newtheorem{problem}{Problem}
\newtheorem*{corollary*}{Corollary}

\title{The Isospectral Problem for $p$-widths: An Application of Zoll Metrics}

\author{Jared Marx-Kuo}
\date{May 2024}
\begin{document}
\maketitle
\begin{abstract}
\noindent We pose the isospectral problem for the $p$-widths: Is a riemannian manifold $(M^n, g)$ uniquely determined by its $p$-widths, $\omega_p(M,g)$? We construct many counterexamples on $S^2$ using Zoll metrics and the fact that geodesic $p$-widths are given by unions of immersed geodesics.
\end{abstract}
\section{Introduction}
In \cite{gromov2006dimension}, Gromov introduced the \textit{$p$-widths}, $\{\omega_p\}$, of a riemannian manifold, $(M^{n+1}, g)$, as an analogue of the spectrum of the laplacian. Intuitively, one replaces the rayleigh quotient in the definition of eigenvalues with the area functional in the definition of $\omega_p$ - see \cite{gromov2010singularities} \cite{gromov2002isoperimetry} \cite{guth2009minimax} \cite{fraser2020applications} for more background. The $p$-widths have proven to be extremely useful due to the Almgren--Pitts/Marques--Neves Morse theory program for the area functional (see \cite{marques2014min} \cite{marques2023morse} \cite{marques2015morse} \cite{marques2021morse} \cite{almgren1962homotopy} \cite{zhou2020multiplicity}). In \cite{liokumovich2018weyl}, Liokumovich--Marqueves--Neves demonstrated that $\{\omega_p\}$ satisfy a Weyl--law, which later lead to the resolution of Yau's conjecture: on any closed manifold $(M^{n+1}, g)$, there exist infinitely many minimal hypersurfaces (see \cite{chodosh2020minimal} \cite{song2018existence} \cite{marques2019equidistribution} \cite{marques2017existence} \cite{irie2018density} \cite{li2023existence}). We refer the reader to the introduction of \cite{chodosh2023p} for a more thorough overview of the history and applications of the p-widths. \nl 
\indent Recall the isospectral problem for the laplacian
\begin{problem} \label{LaplacianIsospectral}
	Is a closed manifold, $(M^n, g)$, uniquely determined by its discrete spectrum, $\{\lambda_k\}$?
\end{problem}
\noindent Informally, this is known as ``can you hear the sound of a drum?" and was made famous by Kac \cite{kac1966can}. In \cite{milnor1964eigenvalues}, Milnor constructed two $16$-dimensional tori which are isospectral but not isometric, providing the first counterexample to \ref{LaplacianIsospectral}. Gordon--Webb--Wolpert \cite{gordon1992isospectral} constructed a simpler counterexample using polygons in $\R^2$. Hersch (\cite{hersch1970quatre}) showed that the round metric on $S^2$ is isospectrally unique by proving an isoperimetric inequality for the first eigenvalue. \nl 
%
%Hersch first eigenvalue reference in here: https://dms.umontreal.ca/~iossif/sphere.pdf
%
\indent Given the analogy between $\{\omega_p\}$ and $\{\lambda_p\}$, 
%\red{Insert introduction, discussion of $p$-widths. Their applications to yau's conjecture, $p$-widths of a surface, etc. Discuss analogy with spectral theory and eigenvalues of the laplacian. Maybe in a later section mention the progress on isospectral problem (not true in general but Steve zelditch showed its true for convex domains)}. \nl \nl
%
%Given the strong analogies and historical motivation between eigenvalues and $\omega_p$, 
It is natural to pose the \underline{\textit{$p$-width isospectral problem}}:
\begin{problem} \label{pWidthIsospectral}
Is a closed manifold, $(M^n, g)$, uniquely determined by the values $\omega_p(g)$?
\end{problem}
\noindent We show that the answer is ``no":
\begin{theorem} \label{constantWidthThm}
For any $p > 0$, $\omega_p$ is constant on any fixed connected component of the space of Zoll metrics on $S^2$.
\end{theorem}
\noindent As a corollary, we can compute the $p$-widths explicitly for new metrics. 
%. Denote $\{(p, \omega_p(g)) \; : \; p \in \Z^+\}$ the \textbf{$p$-width spectrum} of a manifold.
%
\begin{corollary} \label{IsoSpectralThm}
Let $g$ be any Zoll metric on $S^2$ which lies in the connected component of $g_{round}$ in the space of Zoll metrics on $S^2$. Then $\omega_p(g) = \omega_p(g_{round}) = 2 \pi \lfloor \sqrt{p} \rfloor$ for all $p$.
\end{corollary}
%
%\red{Mention how $S^2$ no longer isospectrally rigid! Contrast Li-Yau}
\noindent The above gives rise to a continuous deformation of p-width isospectral metrics (compare \cite{gord1984iso}). In particular, any odd function $f: S^2 \to \R$ gives rise to a $1$-parameter family of geometrically distinct Zoll metrics near the round metric (see \cite{guillemin1976radon}). The $p$-widths of $(S^2, g_{round})$ were computed explicitly in \cite{chodosh2023p}, and this was the first example of the $p$-width spectrum being known for all $p$. Corollary \ref{IsoSpectralThm} now gives uncountably more manifolds where the full width spectrum is known. \nl 
\indent The proof of theorem \ref{constantWidthThm} uses the defining properties of zoll metrics on $S^2$, which are metrics for which all geodesics (with multiplicity one) are simple and length $2 \pi$ (note this choice of length, as opposed to any constant $c$, agrees with \cite{guillemin1976radon}). These were first constructed by Zoll \cite{zoll1901ueber}, and we refer the reader to \cite{guillemin1976radon} \cite{funk1911flachen} \cite{gromoll1981metrics} \cite{besse2012manifolds} for a by no means a complete list of background sources. The proof also utilizes a recent result of Chodosh--Mantoulidis \cite{chodosh2023p}, which states that each $\omega_p(M^2, g)$ is achieved by a union of closed, potentially immersed geodesics. \nl 
\indent This paper is organized as follows: in \S \ref{background}, we briefly define the $p$-widths, and recall the results of Chodosh--Mantoulidis \cite{chodosh2023p}. We also define Zoll metrics and list some relevant properties. In \S \ref{MainThmSec} we prove the main theorem, and in \S \ref{OpenSec} we pose some open questions. The author thanks Otis Chodosh for his conversations on $p$-widths and Lucas Ambrosio for his conversations on Zoll metrics.
%
\begin{comment}
The proof is based on three ideas. One is the existence of zoll metrics on $S^2$, which are metrics such that all geodesics are simple and closed (\red{Check if always simple. Result of gromov and goll}, see \href{https://mathoverflow.net/questions/430669/on-properties-of-besse-spheres}{here}) and have the same length. The second idea is that the $p$-widths on a surface are given by unions of closed geodesics with multiplicity, as proved recently by Chodosh--Mantoulidis (\red{Source}). The final idea is that Zoll metrics are abundant and in fact, uncountable. A counting argument then finishes the proof. 
%
\subsection{Paper Organization}
The paper is organized as follows
%
\begin{itemize}
\item In \S \ref{pwidthBack}, we define the $p$-widths $\omega_p$, and recall the result of Chodosh--Mantoulidis on what types of gedoesics can achieve the $p$-width
\item In \S \ref{zollBack} we define Zoll metrics and recall basic properties of them
\item In \S \ref{MainThmSec} we prove the main theorem
\item In \S \ref{OpenSec} we pose some open questions
\end{itemize}
\end{comment}
%
\section{Background} \label{background}
Let $X \subseteq I^k$ denote a cubical subcomplex. 
\begin{definition}[\cite{marques2017existence}, Defn 4.1] \label{defi:sweepout}
	A map $\Phi : X\to \mathcal{Z}_{1}(M;\Z_{2})$ is a \textbf{$p$-sweepout} if it is continuous (with the standard flat norm topology on $\mathcal{Z}_1(M; \Z_2)$) and $\Phi^{*}(\overline\lambda^{p}) \neq 0$. 
\end{definition}

\begin{definition}[\cite{marques2017existence} \S 3.3] \label{defi:no.concentration.of.mass}
	A map $\Phi : X \to \mathcal{Z}_1(M; \Z_2)$ is said to have \textbf{no concentration of mass} if
	\[ \lim_{r \to 0} \sup \{ \Vert \Phi(x) \Vert(B_r(p)) : x \in X, p \in M \} = 0. \]
\end{definition}

\begin{definition}[\cite{gromov2002isoperimetry}, \cite{guth2009minimax}] \label{defi:p-width}
	We define $\mathcal{P}_{p} = \mathcal{P}_{p}(M)$ to be the set of all $p$-sweepouts, out of any cubical subcomplex $X$, with no concentration of mass. The \textit{$p$-width} of $(M,g)$ is
	\[
	\omega_{p}(M,g) := \inf_{\Phi \in \mathcal{P}_{p}} \sup\{\mathbf{M}(\Phi(x)) : x \in \text{domain}(\Phi)\}.
	\]
\end{definition}
\noindent To see the analogy between $\omega_p$ and $\lambda_p$, consider the following characterization of $\lambda_p(M)$ as
\begin{align*}
	\lambda_p(M) &:= \inf_{\substack{V \subseteq C^{\infty}(M) \\ \dim V = p+1}} \sup_{f \in V \backslash \{0\}} \frac{\int_M |\n f|^2}{\int_M f^2}
\end{align*}
Note the parallels between the min-max definitions of $\omega_p, \lambda_p$, with $\Phi^{*}(\overline\lambda^{p}) \neq 0$ being analogous to $\dim(V) = p+1$. \nl 
\indent Having defined the $p$-widths, note that when $n + 1 = 2$, the widths may be achieved by geodesic networks a priori. In \cite{chodosh2023p}, Chodosh--Mantoulidis used the Allen--Cahn equation to prove that the $p$-widths can always be achieved by a union of closed geodesics:
\begin{theorem}[\cite{chodosh2023p}, Thm 1.2] \label{ChoManPWidths}
	Let $(M^2, g)$ be a closed Riemannian manifold. For every $p \in \Z^+$, there exists a collection of primitive, closed geodesics $\{\sigma_{p,i}\}$ such that 
	\[
	\omega_p(M, g) = \sum_{j = 1}^{N(p)} m_j \cdot L_g(\sigma_{p,j}), \qquad \qquad m_j \in \Z^+
	\]
\end{theorem}
\noindent Here, ``primitive" means that the geodesic is connected and transversed with multiplicity one. They also showed:
\begin{theorem}[\cite{chodosh2023p}, Thm 1.4] \label{SphereWidths}
For $g_0$ the unit round metric on $S^2$
\[
\forall p \in \N, \qquad \omega_p(S^2, g_0) = 2 \pi \lfloor \sqrt{p} \rfloor
\]
\end{theorem}
\noindent We now collect relevant facts about zoll metrics on $S^2$. 
\begin{definition} \label{ZollMetric}
A metric, $g$, on $S^2$	is called a \textit{Zoll metric} if all of the primitive geodesics are closed and length $2 \pi$.
\end{definition}
\noindent Zoll \cite{zoll1901ueber}, Funk \cite{funk1911flachen}, Guillemin \cite{guillemin1976radon}, and Gromoll--Grov \cite{gromoll1981metrics} each contributed significantly to the theory of Zoll metrics. In particular, Gromoll--Grov showed that all of the primitive geodesics on a zoll surface are simple. Furthermore, the combined work of Guillemin and Funk demonstrated that the tangent space of zoll metrics at the round metric is given by odd functions $f: S^2 \to \R$. We also mention a higher dimensional analogue of these metrics constructed by Ambrozio--Marques--Neves \cite{ambrozio2021riemannian}. \nl 
\indent Let $\mathcal{Z}(S^2)$ denote the space of zoll metrics on $S^2$, $\mathcal{Z}_{\alpha}$, an arbitrary connected component of $\mathcal{Z}(S^2)$, and $\mathcal{Z}_0$, the connected component containing $g_0$, the unit round metric. 
\section{Proof of Main theorem} \label{MainThmSec}
Let $\mathcal{M}_{\infty}(S^2)$ denote the set of all smooth metrics on $S^2$, and let $g_0$ denote the round metric. For fixed $p$, $\omega_p: \mathcal{M}_{\infty}(S^2) \to \R$ is continuous due to the following lemma of Marques--Neves--Song:
\begin{lemma}[\cite{marques2019equidistribution}, Lemma 1] \label{WidthContinuity}
Let $C_1 < C_2$. There exists a $K(C_1, C_2)$, such that 
\[
|\omega_p(g) - \omega_p(g')| \leq K(C_1, C_2) p^{1/2} ||g - g'||_{C^0(g_0)}
\]
for any $g, g'$ such that $C_1 g_0 \leq g, g' \leq C_2 g_0$.
\end{lemma}
\noindent Take any zoll metric, $g$, on $S^2$. Applying theorem \ref{ChoManPWidths}, we have
\begin{align*}
	\omega_p(S^2, g) &= \sum_{j = 1}^{N(p)} m_j \cdot L_{g}(\sigma_{p,j}) \\
	&= 2 \pi \sum_{j = 1}^{N(p)} m_j \in 2 \pi \Z^{+}
\end{align*}
since each $\sigma_{p,j}$ is primitive and $L_{g}(\sigma_{p,j}) = 2 \pi$ by definition of Zoll metric. Let $\mathcal{Z}_{\alpha}$ a fixed, connected component of $\mathcal{Z}$. Consider the restriction of $\omega_p: \mathcal{Z}_{\alpha} \to \R^{\infty}$. Its image is discrete, lying inside $2\pi \N$, but $\omega_p$ is continuous in each component due to lemma \ref{WidthContinuity}, thus the map is constant on connected components, proving theorem \ref{constantWidthThm}. 
\section{Open Questions} \label{OpenSec}
Having shown a counterexample to the isospectral $p$-width problem, it would be interesting to know any of the following
\begin{itemize}
%\item Is $(S^2, g_{round})$ rigid for the $p$-width spectrum?
%
%\item We remark that though theorem \ref{IsoSpectralThm} resolves the isospectral problem, we do not explicitly compute the $p$-widths of any of these surfaces. It is natural to find other metrics on $S^2$ which are $p$-width isospectral to the round metric. These were computed in (\red{ref}) as 
%
%\begin{equation} \label{roundSpec}
%\omega_p(g_{round}) = 2 \pi \lfloor \sqrt{p} \rfloor
%\end{equation}
%%
%In particular, if one could show this, then it would demonstrate that the round metric on the sphere is not isospectrally rigid for $p$-widths, in contrast to the behavior of the laplacian.
%
%Local isospectral rigidity?
%
%While $\omega_p$ depends continuously on the choice of metric, the best known rate of convergence is not uniform in $p$ (see \red{ref Song Marques Neves}). That being said, we conjecture that a more clever choice of zoll metrics would give the same spectrum as in \eqref{roundSpec}.
%
\item Do there exist counterexamples to the isospectral problem in dimensions $3 \leq n + 1 \leq 7$?  The $p$-widths behave differently in these dimensions compared to $n+1 = 2$, as demonstrated by theorem \ref{ChoManPWidths}. See \cite{aiex2016width} for further differences relating to multiplicity and index bounds.

\item Do there exist counterexamples to the isospectral problem on other surfaces $(\Sigma^2, g)$?

\item Do there exist other reasonable geometric conditions that one can impose on a manifold so that isospectral rigidity holds?
\end{itemize}
%
%\bibliographystyle{alpha}
%\bibliography{bib}
\bibliography{Isospectral_p_widths_bib}{}
\bibliographystyle{amsalpha}
\end{document}